\documentclass[10pt]{article}

\usepackage[a4paper]{geometry}
\usepackage{indentfirst}
\usepackage[pdfpagemode=UseNone]{hyperref}
\usepackage{amssymb, amsmath, amsthm}

\newcommand{\mynewtheorem}[2]{\newtheorem{#1}{\indent #2}}
\newcommand{\myalttheorem}[2]{\newtheorem*{#1}{\indent #2}}
\mynewtheorem{lemma}{Lemma}
\mynewtheorem{theorem}{Theorem}
\myalttheorem{theorem*}{Theorem}
\newenvironment{myproof}[1][Proof]{\begin{proof}[\indent #1]}{\end{proof}}

\newcommand{\bfa}{\mathbf{a}}
\newcommand{\bfb}{\mathbf{b}}
\newcommand{\bfc}{\mathbf{c}}

\newcommand{\bfe}{\mathbf{e}}

\newcommand{\bfu}{\mathbf{u}}
\newcommand{\bfv}{\mathbf{v}}
\newcommand{\bfw}{\mathbf{w}}
\newcommand{\bfx}{\mathbf{x}}
\newcommand{\bfy}{\mathbf{y}}

\newcommand{\bfo}{\mathbf{0}}

\DeclareMathOperator{\vol}{vol}
\newcommand{\vs}{\vol_\textnormal{Sign}}

\begin{document}

\title{\textbf{Translation Orders on Grids}}
\author{Nikolai Beluhov}
\date{}

\maketitle

\begin{abstract} Let $S$ be a grid in $\mathbb{Z}^d$. We wish to order $S$ so that the way two points compare depends only on the offset between them. Does each such order arise from a linear function? Bhattacharya conjectured, while Ren and Huang confirmed, that in two dimensions the answer is yes. On the other hand, wild orders do exist when $d \ge 4$. We settle the remaining case, by showing that the answer is yes in three dimensions as well. \end{abstract}

\section{Introduction} \label{intro}

Let $S$ be a subset of $\mathbb{Z}^d$. We wish to order $S$ so that the way two points compare depends only on the offset between them. Or, equivalently, the way two points compare must be preserved by translation. Formally, we are looking for a strict total order $\prec$ on $S$ such that $\bfa' \prec \bfb' \Leftrightarrow \bfa'' \prec \bfb''$ whenever $\bfa' + \bfb'' = \bfa'' + \bfb'$. Then we call $\prec$ a \emph{translation order} on $S$.

One simple way to construct translation orders is as follows: Pick any linear function $f : \mathbb{Z}^d \to \mathbb{R}$ injective over $S$, and order the points of $S$ based on their images under $f$. We call such orders \emph{weighted}. Given $S$, it is natural to ask: Do all translation orders on $S$ arise in this manner?

We are interested most of all in the case when $S$ is a grid; i.e., the Cartesian product of finite integer intervals. First, though, we briefly review some other instances of the problem.

The case when $S$ is the Boolean hypercube $\{0, 1\}^d$ (or, equivalently, a grid of size $2 \times 2 \times \cdots \times 2$) has been considered in probability theory. Then \cite{KPS} the answer is ``yes'' when $d \le 4$ but ``no'' otherwise, when $d \ge 5$.

The case when $S$ is the positive orthant $\mathbb{N}^d$, or the entire integer lattice $\mathbb{Z}^d$, has been considered in computer algebra. Then \cite{R} the answer is ``almost'': We cannot always find a suitable $f$, but we can find a collection of linear functions $f_1$, $f_2$, $\ldots$, $f_k$, with $k \le d$, such that the ties of $f_1$ are broken by $f_2$, the ties of $f_2$ are broken by $f_3$, and so on. Notice that, as a corollary, the restriction of a translation order on $\mathbb{N}^d$ or $\mathbb{Z}^d$ to any finite subset of the domain must necessarily be weighted.

The case when $S$ is a discrete simplex of the form $\{(x_1, x_2, \ldots, x_d) \mid x_i \ge 0 \text{ for all } i \text{ and } x_1 + x_2 + \cdots + x_d \le k\}$ has been considered in economics. (Strictly speaking, it is actually the embedding $\{(x_1, x_2, \ldots, x_{d + 1}) \mid x_i \ge 0 \text{ for all } i \text{ and } x_1 + x_2 + \cdots + x_{d + 1} = k\}$ of this discrete simplex in $\mathbb{Z}^{d + 1}$ that has been considered. But the two configurations are affinely isomorphic.) Then \cite{SS} the answer is ``yes'' for all $k$ when $d \le 2$ but ``no'' for some $k$ otherwise, when $d \ge 3$.

For the broader context, we point readers towards the survey \cite{S} and its bibliography.

We focus now on grids. The problem is trivial in one dimension. Then $S$ becomes a finite integer interval, and so it admits exactly two translation orders, both of which are indeed weighted. In two dimensions, Bhattacharya \cite{B} conjectured, while Ren and Huang \cite{RH} confirmed, that the answer is in the affirmative once again. We discuss the two-dimensional setting more carefully in Section \ref{2d}.

We offer a generalisation of the preceding result in Section \ref{ccs}, as follows: We say that $S$ is \emph{centrally symmetric} when it is symmetric across some point of $\mathbb{R}^d$. We also say that $S$ is \emph{convex} when it is the set of all integer points in some convex region of $\mathbb{R}^d$.

\begin{theorem} \label{ccst} Let $S$ be a finite, convex, centrally symmetric set of integer points in the plane. Then every translation order on $S$ is weighted. \end{theorem}

Our main contribution is that we solve the problem for three-dimensional grids as well. The solution occupies Sections \ref{3di}--\ref{3diii}, with the answer still in the affirmative:

\begin{theorem} \label{3dt} {\fontdimen4\font = 1.5pt Let $S$ be a three-dimensional grid. Then every translation order on $S$ is~weighted.} \end{theorem}

There do exist unweighted translation orders on four-dimensional grids, and so a result analogous to Theorem \ref{3dt} cannot hold in four or more dimensions. We discuss the translation orders of higher-dimensional grids more thoroughly in Section \ref{further}.

\section{Initial Observations} \label{init}

We begin with some basic definitions and notations.

A \emph{grid} of size $n_1 \times n_2 \times \cdots \times n_d$ is any subset in $\mathbb{Z}^d$ of the form $I_1 \times I_2 \times \cdots \times I_d$, with each $I_i$ an integer interval of size $n_i$. We call $n_1$, $n_2$, $\ldots$, $n_d$ the \emph{dimensions} of the grid. The \emph{bounding box} of a finite subset in $\mathbb{Z}^d$ is the smallest grid which contains it.

We write $\bfo$ for the origin $(0, 0, \ldots, 0)$ of our coordinate system. The \emph{sign} of a real number is the symbol $-$, $0$, or $+$ when this number is respectively negative, zero, or positive. The \emph{signature} of an ordered tuple of real numbers is the word over the alphabet $\{-, 0, +\}$ formed by their signs.

Let $\bfa_1$, $\bfa_2$, $\ldots$, $\bfa_k$ be affinely independent points in $\mathbb{R}^d$. We write $\bfa_1\bfa_2 \cdots \bfa_k$ for the simplex with these vertices; i.e., the set of all linear combinations $\lambda_1\bfa_1 + \lambda_2\bfa_2 + \cdots + \lambda_k\bfa_k$ with the $\lambda$'s nonnegative and $\lambda_1 + \lambda_2 + \cdots + \lambda_k = 1$. For example, $\bfa_1\bfa_2$ denotes the closed straight-line segment joining $\bfa_1$ and $\bfa_2$.

By the \emph{boundary} of a polytope, we always mean its relative boundary, i.e., its topological boundary within its affine hull; and similarly for the polytope's \emph{interior}. For example, to us the interior of the simplex $\bfa_1\bfa_2 \cdots \bfa_k$ is the set of all linear combinations $\lambda_1\bfa_1 + \lambda_2\bfa_2 + \cdots + \lambda_k\bfa_k$ with the $\lambda$'s positive and $\lambda_1 + \lambda_2 + \cdots + \lambda_k = 1$. We denote the $d$-dimensional volume of a polytope $P$ in $\mathbb{R}^d$ by $\vol P$.

Let $\bfa_1$, $\bfa_2$, $\ldots$, $\bfa_d$ be linearly independent points in $\mathbb{R}^d$. We write $\Pi(\bfa_1, \bfa_2, \ldots, \bfa_d)$ for the parallelotope with one vertex at the origin and its neighbouring vertices at $\bfa_1$, $\bfa_2$, $\ldots$, $\bfa_d$; i.e., the set of all linear combinations $\lambda_1\bfa_1 + \lambda_2\bfa_2 + \cdots + \lambda_d\bfa_d$ with $0 \le \lambda_i \le 1$ for all $i$.

Let $\Pi = \Pi(\bfa_1, \bfa_2, \ldots, \bfa_d)$. Recall that, in this setting, $\vol \bfo\bfa_1\bfa_2 \cdots \bfa_d = 1/d! \cdot \vol \Pi$. We write $\vs \Pi$ for the signed $d$-dimensional volume of $\Pi$, as given by the determinant of the $d \times d$ matrix formed by $\bfa_1$, $\bfa_2$, $\ldots$, $\bfa_d$. (Either as rows, from top to bottom, or as columns, from left to right.) Of course, $\vol \Pi = |\vs \Pi|$.

Let $\bfa_1$, $\bfa_2$, $\ldots$, $\bfa_k$ be arbitrary points in $\mathbb{R}^d$. We write $\operatorname{cone}(\bfa_1, \bfa_2, \ldots, \bfa_k)$ for their conic hull; i.e., the set of all linear combinations $\lambda_1\bfa_1 + \lambda_2\bfa_2 + \cdots + \lambda_k\bfa_k$ with the $\lambda$'s nonnegative.

The plane has its own ``dialect''. Let $\bfa_1$, $\bfa_2$, $\ldots$, $\bfa_k$ be the vertices of a polygon, in this order on its boundary. Then we denote this polygon by $\bfa_1\bfa_2 \cdots \bfa_k$. We also write $\operatorname{area} Q$ for the area of the polygon $Q$.

We go on now to some general observations about translation orders.

Let $S$ be a finite subset of $\mathbb{Z}^d$, with $\prec$ a translation order on $S$. Consider the difference set $D = S - S = \{\bfa - \bfb \mid \bfa, \bfb \in S\}$. We colour each nonzero point $\bfu$ of $D$ black when $\bfa + \bfu = \bfb$ implies $\bfa \prec \bfb$, and white when it implies $\bfa \succ \bfb$. (This sort of construction is standard in similar contexts. See \cite{S}.) What can be said of that colouring?

For each pair of nonzero points $\bfu$ and $-\bfu$ in $D$ symmetric across the origin, one of them must be white and the other one must be black. We call this the \emph{antisymmetry} rule.

Furthermore, if $\bfu_1$, $\bfu_2$, $\bfu_3$ are three nonzero points in $D$ for which there exist three points $\bfa_1$, $\bfa_2$, $\bfa_3$ in $S$ with $\bfa_1 + \bfu_1 = \bfa_2$, $\bfa_2 + \bfu_2 = \bfa_3$, $\bfa_3 + \bfu_3 = \bfa_1$, then $\bfu_1$, $\bfu_2$, $\bfu_3$ cannot be all of the same colour. We call this the \emph{triangle} rule.

These are the comprehensive combinatorial constraints on our colouring, in the sense that a colouring of $D \setminus \{\bfo\}$ in black and white induces a translation order $\prec$ on $S$ if and only if it satisfies both of them.

The triangle rule is subtler, and it merits additional discussion. For a triple $\bfu_1$, $\bfu_2$, $\bfu_3$ of nonzero points in $D$ to fall within its scope, it must hold that $\bfu_1 + \bfu_2 + \bfu_3 = \bfo$. The converse, though, is false in general. We say that $S$ is \emph{triangular} when every triple of nonzero points $\bfu_1$, $\bfu_2$, $\bfu_3$ in $D$ with $\bfu_1 + \bfu_2 + \bfu_3 = \bfo$ admits suitable $\bfa_1$, $\bfa_2$, $\bfa_3$ in $S$.

For a triangular $S$, the triangle rule says simply that three nonzero points in $D$ with sum $\bfo$ can never be all of the same colour. By virtue of the antisymmetry rule, this is equivalent to each colour being closed under addition. Or, in other words, if the sum of two black points in $D$ belongs to $D$, too, then it must be black as well; and similarly for the white points of $D$.

So, in the triangular setting, both of our colouring rules may be formulated entirely in terms of $D$ itself, with no reference to the original set $S$. This means that, if $S$ is triangular, we can safely disregard it and work exclusively within $D$.

We turn next to the matter of weightedness. Clearly, $\prec$ is weighted if and only if there exists a hyperplane through the origin which strictly separates the two colours of $D \setminus \{\bfo\}$. (We can assume without loss of generality that the linear function $f$ in the definition of a weighted order is homogeneous. Then the equation of the separating hyperplane will be $f(\bfx) = 0$.)

By antisymmetry, such a separating hyperplane exists if and only if the origin lies outside of the black points' convex hull. Furthermore, by Carath\'eodory's theorem, the origin lies in the black points' convex hull if and only if there exist $k$ among them, say $\bfu_1$, $\bfu_2$, $\ldots$, $\bfu_k$, with $2 \le k \le d + 1$ and affinely independent, such that the origin lies in the interior of the simplex $\bfu_1\bfu_2 \cdots \bfu_k$. We call any such set of black points \emph{enveloping}. So, in summary, weightedness is equivalent to the non-existence of enveloping sets formed by at most $d + 1$ black points.

We proceed now to relate these observations to some special classes of integer point sets $S$.

If $S$ is a grid, then $D$ will be a grid as well. Explicitly, for a grid of size $n_1 \times n_2 \times \cdots \times n_d$, the difference set becomes the grid $[-(n_1 - 1); n_1 - 1] \times [-(n_2 - 1); n_2 - 1] \times \cdots \times [-(n_d - 1); n_d - 1]$ of size $(2n_1 - 1) \times (2n_2 - 1) \times \cdots \times (2n_d - 1)$. We also get the following:

\begin{lemma} \label{triangle} Every grid is triangular. \end{lemma}

\begin{myproof} Suppose that $S$ is a grid, and let $\bfu_1$, $\bfu_2$, $\bfu_3$ be three nonzero points in $D$ with sum $\bfo$. Take any three points $\bfa_1$, $\bfa_2$, $\bfa_3$ in $\mathbb{Z}^d$ such that $\bfa_1 + \bfu_1 = \bfa_2$, $\bfa_2 + \bfu_2 = \bfa_3$, $\bfa_3 + \bfu_3 = \bfa_1$. For each $1 \le i \le d$, the $i$-th dimension of the bounding box of $\{\bfa_1, \bfa_2, \bfa_3\}$ equals the $i$-th dimension of the bounding box of one among $\{\bfa_1, \bfa_2\}$, $\{\bfa_1, \bfa_3\}$, $\{\bfa_2, \bfa_3\}$, and so it cannot exceed the $i$-th dimension of $S$. \end{myproof}

On the other hand, convexity alone does not imply triangularity, and neither does its conjunction with central symmetry. One counterexample is given by $S = \{\bfo, \pm (1, 0), \pm (0, 1), \pm (1, 1)\}$ and the triple $(2, 0)$, $(0, 2)$, $(-2, -2)$ of nonzero points in $D$. Still, there are some useful things we can say already in this generality:

\begin{lemma} \label{pair} There are no enveloping pairs in $D$ when $S$ is convex. \end{lemma}

\begin{myproof} Suppose, for the sake of contradiction, that $\bfu_1$ and $\bfu_2$ form an enveloping pair of black points in $D$. Then $\lambda_1\bfu_1 + \lambda_2\bfu_2 = \bfo$, with $\lambda_1$ and $\lambda_2$ being relatively prime positive integers. Let $\bfu = 1/\lambda_2 \cdot \bfu_1 = -1/\lambda_1 \cdot \bfu_2$.

Since $\bfu_1$ is in $D$, there are $\bfa$ and $\bfb$ in $S$ with $\bfa + \bfu_1 = \bfb$. Since $S$ is convex, all of the intermediate points $\bfa + \bfu$, $\bfa + 2\bfu$, $\ldots$, $\bfa + (\lambda_2 - 1)\bfu$ are in $S$ as well. Since $\prec$ is a translation order, it follows that the colour of $\bfu$ matches the colour of $\bfu_1$. Similarly, the colour of $-\bfu$ matches the colour of $\bfu_2$. However, this makes both of $\bfu$ and $-\bfu$ black, in contradiction with antisymmetry. \end{myproof}

We return to this train of thought in Section \ref{3di}.

\section{Two-Dimensional Grids} \label{2d}

Here, we briefly review the earlier work on two-dimensional grids done in \cite{B} and \cite{RH}. This work is phrased in terms of a different but isomorphic family of combinatorial objects, and so we begin by mapping out the isomorphism.

Fix two positive integers $m$ and $n$. Consider a grid $G$ of size $(2m - 1) \times (2n - 1)$ whose points are coloured in black and white. There are a total of $mn$ subgrids in $G$ of size $m \times n$. We call each one of them a \emph{window}, and for each window we count the black points in it.

Bhattacharya \cite{B} calls a colouring of $G$ \emph{diverse} when these counts are pairwise distinct, and conjectures that the diverse colourings of $G$ coincide with the colourings of $G$ obtained as follows: Draw a line $\ell$ through the centre of $G$ which does not meet $G$ elsewhere, colour this centre arbitrarily, and colour the other points of $G$ based on which half-plane relative to $\ell$ contains them.

We claim that this proposition is equivalent to the following: The translation orders on the grid $H$ of size $m \times n$ are all weighted.

Indeed, observe first that the centre of $G$ is contained in all windows, and so nothing changes if we make it grey. There are now $mn$ windows and $mn$ options for number of black points in a window, namely $0$, $1$, $\ldots$, $mn - 1$. This means that each option is realised exactly once. Similarly, for each one among $0$, $1$, $\ldots$, $mn - 1$, exactly one window contains this many white points.

We can assume without loss of generality that $G$ is the difference set of $H$, i.e., $[-(m - 1);\allowbreak {m - 1}] \times [-(n - 1); n - 1]$. Consider next the grid $U = [1; mn] \times [1; mn]$ of size $mn \times mn$. We lift the colouring of $G$ into a colouring of $U$, as follows: Fix any enumeration $\bfa_1$, $\bfa_2$, $\ldots$, $\bfa_{mn}$ of $H$. We colour each point $(i, j)$ of $U$ in the same colour as the point $\bfa_i - \bfa_j$ of $G$. So, in particular, all points $(i, i)$ on the main diagonal of $G$ become grey. Observe also that each column $i$ of $U$ is coloured as the window $\bfa_i - H = \{\bfa_i - \bfa \mid \bfa \in H\}$ in $G$; and, similarly, each row $j$ of $U$ is coloured as the window $H - \bfa_j = \{\bfa - \bfa_j \mid \bfa \in H\}$ in $G$.

The column of $U$ with no black points and the row of $U$ with no white points must meet on grey. Following this, the column of $U$ with a single black point and the row of $U$ with a single white point must meet on grey as well. Continuing on in this way, we get that the column of $U$ with $i$ black points and the row of $U$ with $i$ white points must meet on grey for all $i$. By revising the original enumeration of $H$ as needed, we can ensure for all $i$ that these are, in fact, column $i + 1$ and row $i + 1$ of $U$. But this makes each nonzero point $\bfa_i - \bfa_j$ of $G$ white when $i < j$ and black when $i > j$. Hence, our enumeration of $H$ is in fact a translation order, with the colouring of $G$ being precisely the colouring associated with this translation order as in Section~\ref{init}.

The aforementioned conjecture was confirmed by Ren and Huang \cite{RH}. We proceed now to sketch their proof, in broad strokes. By the preceding discussion, this proof also confirms the weightedness of translation orders on two-dimensional grids.

The points of $G$ are coloured in batches, step by step. The intermediate states of $G$ during this process are called \emph{partial colourings} -- with some points white, some black, and the rest uncoloured. Each one of these partial colourings is characterised by a certain ordered triple of vectors with nonnegative integer components. Let $\Psi(\bfu_1, \bfu_2, \bfu_3)$ denote the partial colouring characterised by $\bfu_1$, $\bfu_2$, $\bfu_3 = (x_3, y_3)$. In it, all points $(x, y)$ of $G$ with $x \ge x_3$ and $y \ge y_3$ are left uncoloured; so are their reflections across the origin; of the other points of $G$, each one $\bfv$ with $\operatorname{area}_\text{Sign} \Pi(\bfu_1 + \bfu_2, \bfv) > 0$ is made white; and each one with $\operatorname{area}_\text{Sign} \Pi(\bfu_1 + \bfu_2, \bfv) < 0$ is made black.

Let $\bfe_1 = (1, 0)$ and $\bfe_2 = (0, 1)$ be the basis vectors of our coordinate system. The initial partial colouring that the proof begins with is $\Psi(\bfe_1, \bfe_2, \bfe_1 + \bfe_2)$. Then a sequence of nested partial colourings is developed out of it, with the help of the following claims: (i) Suppose that $\operatorname{area}_\text{Sign} \Pi(\bfu_1 + \bfu_2, \bfu_3) > 0$. Then any completion of $\Psi(\bfu_1, \bfu_2, \bfu_3)$ is also a completion of $\Psi(\bfu_1, \bfu_2, \bfu_3 + \bfe_1)$. (ii) Suppose that $\operatorname{area}_\text{Sign} \Pi(\bfu_1 + \bfu_2, \bfu_3) < 0$. Then any completion of $\Psi(\bfu_1, \bfu_2, \bfu_3)$ is also a completion of $\Psi(\bfu_1, \bfu_2, \bfu_3 + \bfe_2)$. (iii) Suppose that $\bfu_1 + \bfu_2 = \bfu_3$. Then any completion of $\Psi(\bfu_1, \bfu_2, \bfu_3)$ is also a completion of either $\Psi(\bfu_1, \bfu_3, \bfu_3 + \bfe_1)$ or $\Psi(\bfu_3, \bfu_2, \bfu_3 + \bfe_2)$.

These claims are established by considering the following question: Which window in $G$ is going to have as many black points in the complete colouring as $H' = [x_3 - (m - 1); x_3] \times [0; n - 1]$ or $H'' = [-(m - 1); 0] \times [-y_3; -y_3 + (n - 1)]$ does in the partial colouring? By applying claims (i), (ii), (iii) iteratively until no uncoloured points remain, in the end all diverse colourings are demonstrated to be as desired. This completes the sketch of Ren and Huang's proof.

Of course, the weightedness of translation orders on two-dimensional grids is also a special case of both Theorems \ref{ccst} and \ref{3dt}. We present an additional direct proof of it in Section \ref{3di}.

\section{Convexity and Central Symmetry} \label{ccs}

Here, we prove Theorem \ref{ccst}.

Note that neither the central symmetry not the convexity may be waived. For example, the set $\{(-2, 0), (-1, 0), (1, 0), (2, 0)\}$ is centrally symmetric but not convex; the set $\{(0, 0), (1, 0),\allowbreak (0, 1), (-1, -1)\}$ is convex but not centrally symmetric; and both of them admit unweighted translation orders.

Furthermore, the result does not generalise to higher dimensions. For example, the set of all integer points in $\Pi((1, 0, 1), (0, 1, 1), (-1, -1, 1))$ admits the unweighted translation order $(0, 0, 0) \prec (0, 0, 1) \prec (1, 0, 1) \prec (-1, -1, 1) \prec (0, -1, 2) \prec (0, 1, 1) \prec (1, 1, 2) \prec (-1, 0, 2) \prec (0, 0, 2) \prec (0, 0, 3)$. So, in particular, we cannot salvage the higher-dimensional analogues by restriction to the setting where the convex hull of $S$ is a zonotope or a parallelotope.

We go on now to the proof. It proceeds by induction on the size of $S$.

Observe, to begin with, that if the points of $S$ are collinear, then $S$ becomes essentially one-dimensional and the result is clear. This deals away with the base cases of our induction, where $|S| \le 3$.

For the induction step, we assume that $|S| \ge 4$ and the convex hull $Q$ of $S$ is properly two-dimensional. So $Q$ is a centrally symmetric convex polygon with integer vertices and $S$ is the set of all integer points in $Q$. We say that two points of $S$ are \emph{opposites} when they are symmetric across the centre $\bfc$ of $Q$.

\begin{lemma} \label{vr} Let $\bfa$ be a point of $S$. Then there is some vertex $\bfu$ of $Q$ such that the reflection of $\bfu$ across $\bfa$ belongs to $S$. \end{lemma} 

\begin{myproof} Any vertex of $Q$ works when $\bfa = \bfc$. Otherwise, the line through $\bfc$ and $\bfa$ pierces two opposite sides of $Q$, say $\bfv\bfw$ and $\bfv'\bfw'$. Denote the parallelogram $\bfv\bfw\bfv'\bfw'$ by $R$.

The two lines through $\bfc$ parallel to the sides of $R$ partition it into four smaller parallelograms. Consider the one which contains $\bfa$, and let $\bfu$ be the unique vertex of $R$ in it. Then the reflection of $\bfu$ across $\bfa$ will be an integer point in $R$, and hence also an element of $S$. \end{myproof}

\begin{lemma} \label{do} The $\prec$-smallest and $\prec$-greatest elements of $S$ are opposite vertices of $Q$. \end{lemma} 

\begin{myproof} Let $\bfa$ be the $\prec$-smallest element of $S$. Suppose, for the sake of contradiction, that it is not a vertex of $Q$. By Lemma \ref{vr}, there is some vertex $\bfu$ of $Q$ such that the reflection $\bfv$ of $\bfu$ across $\bfa$ is also an element of $S$. But now $\bfa \prec \bfu$ and $\bfa \prec \bfv$ while $\bfu + \bfv = 2\bfa$, in contradiction with $\prec$ being a translation order. So $\bfa$ must indeed be a vertex of $Q$.

Similarly, the $\prec$-greatest element $\bfb$ of $S$ must be a vertex of $Q$, too. Let $\bfa'$ and $\bfb'$ be the opposites of $\bfa$ and $\bfb$, respectively. Suppose, for the sake of contradiction, that $\bfa' \neq \bfb$ and $\bfb' \neq \bfa$. Then $\bfa \prec \bfb'$ and $\bfa' \prec \bfb$ while $\bfa + \bfa' = \bfb + \bfb'$, in contradiction with $\prec$ being a translation order. So $\bfa$ and $\bfb$ must indeed be opposites. \end{myproof}

For the sequel, it will be convenient to assume that the $\prec$-smallest element of $S$ is the origin~$\bfo$. To maintain notational consistency, we also write $\bfo'$ for the $\prec$-greatest element of $S$. Let $S^\star = S \setminus \{\bfo, \bfo'\}$. By the induction hypothesis, Theorem \ref{ccst} holds for the restriction of $\prec$ to~$S^\star$. Our goal is to put $\bfo$ and $\bfo'$ back in.

The rest of the proof will be easier to express in terms of the black-and-white colouring of $D$ rather than the original translation order on $S$. Notice that, since the origin is the $\prec$-smallest element of $S$, all nonzero points of $S$ are black in $D$ and all nonzero points of $-S = \{-\bfa \mid \bfa \in S\}$ are white in $D$.

Let $D^\star$ be the difference set of $S^\star$. So $D^\star \subseteq D$, and the colour of each point in $D^\star$ coincides with its colour in $D$. We also know that each point in $D \setminus D^\star$ is of the form either $\pm(\bfo - \bfa)$ or $\pm(\bfo' - \bfa)$, with $\bfa \in S$. However, $\pm(\bfo' - \bfa) = \mp(\bfo - \bfa')$, with $\bfa'$ being the opposite of $\bfa$. We conclude that all black points in $D \setminus D^\star$ belong to $S$ and all white points in $D \setminus D^\star$ belong to~$-S$.

Let $\bfv$ be the integer point on the boundary of $Q$ immediately before $\bfo$, going counterclockwise, and similarly let $\bfw$ be the integer point on the boundary of $Q$ immediately after $\bfo$. (For each one of $\bfv$ and $\bfw$, notice that it is not necessarily a vertex of $Q$.) Then the opposites $\bfv'$ and $\bfw'$ of $\bfv$ and $\bfw$ will be, respectively, the integer point immediately before and the integer point immediately after $\bfo'$ on the boundary of $Q$, going counterclockwise once again.

Observe that it suffices to find a separator $\ell$ for $D^\star$ such that both of $\bfv$ and $\bfw$ are in the black half-plane relative to $\ell$. Indeed, since $S \subseteq \operatorname{cone}(\bfv, \bfw)$, it will then follow that $S \setminus \{\bfo\}$ is also contained within the black half-plane relative to $\ell$, and similarly $-S \setminus \{\bfo\}$ will be contained within the white half-plane relative to $\ell$. So $\ell$ will be a separator for $D$ as well.

We consider two cases, based on whether $\bfv$ and $\bfw$ are opposites or not. They are handled by Lemmas \ref{hex}~and~\ref{quad} below, respectively.

\begin{lemma} \label{hex} Suppose that $\bfv$ and $\bfw$ are not opposites. Then every separator for $D^\star$ is also a separator for $D$. \end{lemma}

\begin{myproof} In this case, $\bfv\bfo\bfw\bfv'\bfo'\bfw'$ is a centrally symmetric convex hexagon. So the integer point $\bfv + \bfw$ lies in its interior, and hence belongs to $S^\star$. Since $\bfv$ and $\bfw$ do, too, we get that both of $\bfv = (\bfv + \bfw) - \bfw$ and $\bfw = (\bfv + \bfw) - \bfv$ are elements of $D^\star$. Furthermore, both of them must be black, due to being elements of $S$. By the preceding discussion, this suffices. \end{myproof}

For the other case, we are going to need one auxiliary observation:

\begin{lemma} \label{qa} Let $\bfa$ and $\bfb$ be two linearly independent integer points in the plane, with $R = \Pi(\bfa, \bfb)$. Suppose that there is no integer point on the boundary of $R$ other than its vertices, but in the interior of $R$ there is some integer point outside of the diagonal $\bfa\bfb$. Then we can find two integer points in the interior of $R$ whose difference also belongs to the interior of $R$. \end{lemma} 

\begin{myproof} The two lines through the centre of $R$ parallel to its sides partition it into four smaller parallelograms. For each vertex $\bfu$ of $R$, let $R_\bfu$ be the smaller parallelogram which contains $\bfu$.

Consider the integer point $\bfx$ in the interior of $R$ furthest away from the line through $\bfa$~and~$\bfb$. Suppose, for the sake of contradiction, that it lies in the interior of some $R_\bfu$ with $\bfu \in \{\bfa, \bfb\}$. Then the reflection of $\bfu$ across $\bfx$ will be an integer point in the interior of $R$, too, and it will lie twice as far away from the line through $\bfa$ and $\bfb$, in contradiction with our choice of $\bfx$. So $\bfx$ must lie in some $R_\bfu$ with $\bfu \in \{\bfo, \bfa + \bfb\}$.

Notice that $\bfx$ cannot lie on the boundary of its $R_\bfu$, as otherwise the reflection of $\bfu$ across $\bfx$ will be an integer point on the boundary of $R$ distinct from its vertices. Suppose, without loss of generality, that $\bfx$ lies in the interior of $R_\bfo$ and its reflection $\bfy$ across the centre of $R$ lies in the interior of $R_{\bfa + \bfb}$. Then the difference $\bfy - \bfx$ will belong to the interior of $R$, as needed. \end{myproof}

\begin{lemma} \label{quad} Suppose that $\bfv$ and $\bfw$ are opposites. Then some separator for $D^\star$ is also a separator for $D$. \end{lemma}

\begin{myproof} In this case, $Q$ is the parallelogram $\bfv\bfo\bfw\bfo'$ and there are no integer points on its boundary other than its vertices. So $D \cap \operatorname{cone}(\bfv, \bfw) = S$ and all nonzero points of $D^\star$ in $\operatorname{cone}(\bfv, \bfw)$ belong to the interior of $Q$. Since these points are all black, due to being elements of $S$, we get that there exists some separator $\ell$ for $D^\star$ which does not meet $\operatorname{cone}(\bfv, \bfw)$ except at the origin.

We are done if $D^\star \cap \operatorname{cone}(\bfv, \bfw) \neq \varnothing$, as in this scenario $\operatorname{cone}(\bfv, \bfw)$ must be contained within the black half-plane relative to $\ell$ and the desired conclusion follows as previously discussed. Otherwise, by Lemma \ref{qa} all integer points in the interior of $Q$ must lie on the segment $\bfv\bfw$, implying that $S^\star$ is a subset of that segment while $D^\star$ is a subset of the line through the origin parallel to it. But in that scenario all other lines through the origin become separators for $D^\star$, and so we can easily find one among them whose black half-plane contains $\operatorname{cone}(\bfv, \bfw)$. \end{myproof}

Together, Lemmas \ref{hex} and \ref{quad} complete the induction step as well as the proof of Theorem~\ref{ccst}.

\section{Three-Dimensional Grids I} \label{3di}

Here, we begin to establish Theorem \ref{3dt}.

Let $S$ be a grid. By Lemma \ref{triangle}, $S$ is triangular, and so we can make use of the two simplified forms of the triangle rule as per Section \ref{init}. Since $S$ is convex, we also get that there are no enveloping pairs in $D$, by Lemma \ref{pair}.

Our strategy for the proof of Theorem \ref{3dt} will be to rule out the existence of enveloping triples and quadruples as well. We do this by assuming the converse, for the sake of contradiction, and then finding a way to ``shrink'' any enveloping set. There are two distinct kinds of shrinking that we pursue. One is to find an enveloping set with fewer points. The other is to find an enveloping set with the same number of points, but whose convex hull is of smaller volume. (For the appropriate notion of volume; i.e., the area for triples and the three-dimensional volume for quadruples.)

Suppose we can show that every enveloping triple or quadruple is shrinkable. Then we are done: We just take the enveloping set of smallest size, and with a convex hull of the smallest volume among all enveloping sets of that size. By shrinking it, we arrive at the desired contradiction.

We proceed now to deal away with the triples:

\begin{lemma} \label{triple} There are no enveloping triples in $D$ when $S$ is a grid. \end{lemma}

\begin{myproof} Suppose that the black points $\bfu_1$, $\bfu_2$, $\bfu_3$ form an enveloping triple in $D$. We aim to show that this triple is shrinkable.

Let $\lambda_1$, $\lambda_2$, $\lambda_3$ be positive real numbers with $\lambda_1\bfu_1 + \lambda_2\bfu_2 + \lambda_3\bfu_3 = \bfo$. We can assume, without loss of generality, that $\lambda_1 = \min\{\lambda_1, \lambda_2, \lambda_3\}$. Let $\bfv = \bfu_2 + \bfu_3$. We claim that $\bfv$ belongs~to~$D$.

Indeed, take any $i$ with $1 \le i \le d$ and consider the $i$-th components $u_1$, $u_2$, $u_3$, $v$ of $\bfu_1$, $\bfu_2$, $\bfu_3$, $\bfv$. Of course, $v = u_2 + u_3$. We branch into two cases next, based on the signs of $u_2$ and $u_3$. If these signs differ, $|v| \le \max\{|u_2|, |u_3|\}$. Otherwise, they coincide. We can assume, without loss of generality, that both of $u_2$ and $u_3$ are nonnegative. Then $v = u_2 + u_3 \le \lambda_2/\lambda_1 \cdot u_2 + \lambda_3/\lambda_1 \cdot u_3 = -u_1$. Either way, $|v| \le \max\{|u_1|, |u_2|, |u_3|\}$. Since our reasoning applies to all coordinate positions $i$, it follows that $\bfv$ belongs to the bounding box of the six points $\pm \bfu_1$, $\pm \bfu_2$, $\pm \bfu_3$, and hence also to $D$.

By the triangle rule, $\bfv$ must be black. We claim, finally, that $\bfv$ allows us to shrink the enveloping triple formed by $\bfu_1$, $\bfu_2$, $\bfu_3$.

This is clear when $\lambda_2 = \lambda_3$, as in that case $\bfu_1$ and $\bfv$ form an enveloping pair. Otherwise, we can assume without loss of generality that $\lambda_2 > \lambda_3$. Consider the triple $\bfu_1$, $\bfu_2$, $\bfv$. It is enveloping as $\lambda_1\bfu_1 + (\lambda_2 - \lambda_3)\bfu_2 + \lambda_3\bfv = \bfo$. Furthermore, $\operatorname{area} \bfu_1\bfu_2\bfo + \operatorname{area} \bfu_1\bfu_2\bfv = \operatorname{area} \bfu_1\bfu_2\bfu_2 + \operatorname{area} \bfu_1\bfu_2\bfu_3$ by the linearity of the signed area, and so $\operatorname{area} \bfu_1\bfu_2\bfv < \operatorname{area} \bfu_1\bfu_2\bfu_3$. This confirms the desired shrinking. \end{myproof}

Notice that Lemmas \ref{pair} and \ref{triple} together resolve the problem for two-dimensional grids, as in two dimensions it is only the pairs and triples that matter.

Suppose, from now on, that $S$ is a three-dimensional grid specifically. For the proof of Theorem \ref{3dt}, it remains to rule out the existence of enveloping quadruples. Let $\bfu_1$, $\bfu_2$, $\bfu_3$, $\bfu_4$ be any enveloping quadruple of black points in $D$. Our goal is to shrink it. We denote the tetrahedron $\bfu_1\bfu_2\bfu_3\bfu_4$ by $T$.

\begin{lemma} \label{rep} Let $\bfv$ be any point distinct from $\bfo$, $\bfu_1$, $\bfu_2$, $\bfu_3$, $\bfu_4$. Then we can uniquely discard one or more elements of the set $\{\bfu_1, \bfu_2, \bfu_3, \bfu_4\}$ and replace them with $\bfv$ so that the origin continues to lie in the interior of the convex hull of the new set obtained after the replacement. \end{lemma}

\begin{myproof} Consider the ray emanating from the origin opposite $\bfv$; i.e., the set of all points $\lambda\bfv$ with $\lambda \le 0$. Since the origin lies in the interior of $T$, this ray pierces the boundary of $T$ at a unique point. If this point is a vertex $\bfu'$ of $T$, our new set will be $\{\bfu', \bfv\}$. If it lies in the interior of some edge $\bfu'\bfu''$ of $T$, our new set will be $\{\bfu', \bfu'', \bfv\}$. Finally, if it lies in the interior of some face $\bfu'\bfu''\bfu'''$ of $T$, our new set will be $\{\bfu', \bfu'', \bfu''', \bfv\}$. \end{myproof}

We call the new set thus obtained the \emph{replacement set} of $\bfv$, and its convex hull the \emph{replacement figure} of $\bfv$. We say that $\bfv$ is \emph{good} when it is distinct from $\bfo$, $\bfu_1$, $\bfu_2$, $\bfu_3$, $\bfu_4$ and either its replacement set contains at most three points or the volume of its replacement figure is smaller than the volume of $T$. We say that $\bfv$ is \emph{excellent} when both of $\bfv$ and $-\bfv$ are good.

Let $B$ be the bounding box of the eight points $\pm \bfu_1$, $\pm \bfu_2$, $\pm \bfu_3$, $\pm \bfu_4$. Observe that it suffices to find an excellent point in $B$. Indeed, since $B \subseteq D$, any such excellent point will belong to $D$ as well. By antisymmetry, either it or its reflection across the origin will be black. So we will certainly get an opportunity to shrink our enveloping quadruple, one way or the other.

We turn now to a search for an excellent point in $B$. (Notice that, from here on out, we can safely disregard $D$ and its black-and-white colouring.)

Let $s_1$, $s_2$, $s_3$, $s_4$ be the volumes of the four tetrahedra $\bfo\bfu_2\bfu_3\bfu_4$, $\bfo\bfu_1\bfu_3\bfu_4$, $\bfo\bfu_1\bfu_2\bfu_4$, $\bfo\bfu_1\bfu_2\bfu_3$, respectively. Then the volume $t$ of $T$ is given by $s_1 + s_2 + s_3 + s_4$. Furthermore, since $s_1$, $s_2$, $s_3$, $s_4$ are proportional to the barycentric coordinates of the origin relative to $T$, we get that $s_1\bfu_1 + s_2\bfu_2 + s_3\bfu_3 + s_4\bfu_4 = \bfo$.

We can assume without loss of generality that $s_4 = \max\{s_1, s_2, s_3, s_4\}$. Let $P$ be the parallelepiped $\Pi(\bfu_1, \bfu_2, \bfu_3)$. Then the volume $p$ of $P$ is given by $6s_4$.

For convenience, we sometimes write $\bfu_{12}$ instead of $\bfu_1 + \bfu_2$, $\bfu_{123}$ instead of $\bfu_1 + \bfu_2 + \bfu_3$, and so on. For notational consistency, we might also sometimes write $\bfu_\varnothing$ instead of $\bfo$. Using this shorthand notation, the vertices of $P$ become $\bfu_\varnothing$, $\bfu_1$, $\bfu_2$, $\bfu_3$, $\bfu_{12}$, $\bfu_{13}$, $\bfu_{23}$, $\bfu_{123}$. We write $\iota$ for a generic subscript out of these eight, so that $\bfu_\iota$ stands for a generic vertex of $P$.

We consider, alongside our main coordinate system, also an alternative coordinate system $\Lambda$ of the same origin but with basis vectors $\bfu_1$, $\bfu_2$, $\bfu_3$. We write $\langle \lambda_1, \lambda_2, \lambda_3 \rangle$ for the point $\lambda_1\bfu_1 + \lambda_2\bfu_2 + \lambda_3\bfu_3$; i.e., the point with these coordinates relative to $\Lambda$. So, in particular, $P$ is the set of all points $\langle \lambda_1, \lambda_2, \lambda_3 \rangle$ with $0 \le \lambda_i \le 1$ for all $i$, from the perspective of $\Lambda$. Furthermore, $\bfu_4 = \langle -s_1/s_4, -s_2/s_4, -s_3/s_4 \rangle$, implying that $-\bfu_4$ belongs to $P$.

\begin{lemma} \label{surf} Suppose that there is an integer point on the boundary of $P$ other than its vertices. Then we can find an excellent point in $B$. \end{lemma}

\begin{myproof} Let $\bfv$ be an integer point on the boundary of $P$ other than its vertices. We can assume without loss of generality that $\bfv$ lies on one of the three faces of $P$ which meet at the origin; otherwise, we reflect $\bfv$ across the centre of $P$. We can next assume without loss of generality that $\bfv$ lies on the face $\bfo\bfu_1\bfu_{12}\bfu_2$ of $P$ specifically. Finally, we can also assume without loss of generality that $\bfv$ lies in the triangle $\bfo\bfu_1\bfu_2$; otherwise, we reflect $\bfv$ across the centre of the parallelogram $\bfo\bfu_1\bfu_{12}\bfu_2$.

Since all three of $\bfo$, $\bfu_1$, $\bfu_2$ belong to $B$, and $\bfv$ lies in their convex hull, we get that $\bfv$ belongs to $B$ as well. We proceed to show that it is excellent.

Clearly, $-\bfv$ is good as the origin lies in the convex hull of $\bfu_1$, $\bfu_2$, $-\bfv$, and so the replacement set of $-\bfv$ is either a pair or a triple.

What about $\bfv$? It will be good, too, if its replacement set is either a pair or a triple. Suppose, otherwise, that it is a quadruple. This quadruple cannot contain both of $\bfu_1$ and $\bfu_2$ simultaneously. So the replacement figure of $\bfv$ must be one of the two tetrahedra $\bfu_1\bfu_3\bfu_4\bfv$ and $\bfu_2\bfu_3\bfu_4\bfv$.

Since $\bfv$ lies in the convex hull of $\bfo$, $\bfu_1$, $\bfu_2$, we get that $\bfv = \lambda_1\bfu_1 + \lambda_2\bfu_2 + \lambda_3\bfo$ with $0 \le \lambda_i \le 1$ for all $i$ and $\lambda_1 + \lambda_2 + \lambda_3 = 1$. By the linearity of the signed volume, $\vol \bfu_1\bfu_3\bfu_4\bfv = \lambda_1 \vol \bfu_1\bfu_3\bfu_4\bfu_1 + \lambda_2 \vol \bfu_1\bfu_3\bfu_4\bfu_2 + \lambda_3 \vol \bfu_1\bfu_3\bfu_4\bfo = \lambda_2t + \lambda_3s_2 < t$, as $s_2 < t$ and $\lambda_2 + \lambda_3 \le 1$. Similarly, $\vol \bfu_2\bfu_3\bfu_4\bfv < t$ as well, and the shrinking has been guaranteed. \end{myproof}

In light of Lemma \ref{surf}, we can assume from now on that there are no integer points on the boundary of $P$ other than its vertices.

Consider next the integer points in the interior of $P$. Our search for an excellent point in $B$ will now branch into two rather different cases based on how the set of these integer points is structured. We ask: Do all integer points in the interior of $P$ lie on its diagonal $\bfo\bfu_{123}$, or not? The two scenarios will be resolved in Sections \ref{3dii} and \ref{3diii}, respectively.

\section{Three-Dimensional Grids II} \label{3dii}

Suppose, throughout this section, that there is an integer point in the interior of $P$ outside of its diagonal $\bfo\bfu_{123}$, and also that there are no integer points on the boundary of $P$ other than its vertices.

Let $O$ be the convex hull of the six points $\pm \bfu_1$, $\pm \bfu_2$, $\pm \bfu_3$. So, a convex polyhedron of the same combinatorial type as the regular octahedron. (Indeed, affinely isomorphic to the regular octahedron.) Or, from the perspective of $\Lambda$, the set of all points $\langle \lambda_1, \lambda_2, \lambda_3 \rangle$ with $|\lambda_1| + |\lambda_2| + |\lambda_3|~\le~1$. Clearly, all of its integer points belong to $B$. We are going to search for our excellent point inside of $O$.

\begin{lemma} \label{oct} Let $\bfv$ be any point of $O$ and let $F$ be any face of $T$ other than $\bfu_1\bfu_2\bfu_3$. Then $\vol \bfv F \le \vol T$, with equality attained if and only if $\bfv$ is the vertex of $T$ opposite $F$. \end{lemma}

\begin{myproof} Suppose, for concreteness, that $F = \bfu_2\bfu_3\bfu_4$ and let $f(\bfv) = \vs \bfv\bfu_2\bfu_3\bfu_4 = 1/6 \cdot \vs \Pi(\bfu_2 - \bfv, \bfu_3 - \bfv, \bfu_4 - \bfv)$. By linearity, the extreme values of $f$ as $\bfv$ varies over $O$ must be attained at the vertices of $O$. Let us analyse these vertices.

We can assume without loss of generality that $f(\bfo) = s_1$ and $f(\bfu_1) = t$. By the linearity of the signed volume, $f(\bfu_1) + f(-\bfu_1) = 2f(\bfo)$ and so $f(-\bfu_1) = 2s_1 - t$. Since $0 < s_1 < t$, it follows that $|f(-\bfu_1)| < t$.

Of course, $f(\bfu_2) = 0$. By the linearity of the signed volume, $f(\bfu_2) + f(-\bfu_2) = 2f(\bfo)$ and so $f(-\bfu_2) = 2s_1$. Since $0 < s_1 \le s_4$, it follows that $0 < f(-\bfu_2) < s_1 + s_2 + s_3 + s_4 = t$. Similarly, $f(\bfu_3) = 0$ and $0 < f(-\bfu_3) < t$ as well. \end{myproof}

We now exclude ``two eighths'' out of $O$. Let $H$ be the set of all points in $O$ such that the point's signature relative to $\Lambda$ is distinct from ${-}{-}{-}$ and ${+}{+}{+}$. Or, equivalently: Writing $\bfu'_i$ for $-\bfu_i$, we may describe $H$ as the union of the six tetrahedra $\bfo\bfu'_1\bfu_2\bfu_3$, $\bfo\bfu_1\bfu'_2\bfu_3$, $\bfo\bfu_1\bfu_2\bfu'_3$, $\bfo\bfu'_1\bfu'_2\bfu_3$, $\bfo\bfu'_1\bfu_2\bfu'_3$, $\bfo\bfu_1\bfu'_2\bfu'_3$.

\begin{lemma} \label{68} Every integer point in $H$ distinct from $\bfo$, $\pm \bfu_1$, $\pm \bfu_2$, $\pm \bfu_3$ is excellent. \end{lemma}

\begin{myproof} Let $\bfv$ be any integer point as described. Then $\bfv$ is distinct from $\pm \bfu_4$, too, as the signatures of these points relative to $\Lambda$ put them outside of $H$. Since $H$ is symmetric across the origin, it suffices to show that $\bfv$ is good.

This is automatic when the replacement set of $\bfv$ is either a pair or a triple. Suppose, otherwise, that it is a quadruple, and so the replacement figure of $\bfv$ is a tetrahedron. By the definition of $H$ and the proof of Lemma \ref{rep}, this tetrahedron cannot be $\bfu_1\bfu_2\bfu_3\bfv$. So it must be of the form $\bfv F$, for some face $F$ of $T$ distinct from $\bfu_1\bfu_2\bfu_3$. The desired conclusion now follows by Lemma~\ref{oct}. \end{myproof}

We are only left to find an integer point in $H$ distinct from $\bfo$, $\pm \bfu_1$, $\pm \bfu_2$, $\pm \bfu_3$.

For each vertex $\bfu_\iota$ of $P$, let $\tau_\iota$ be the homothety with centre $\bfu_\iota$ and scaling factor $2$. Or, in other words, $\tau_\iota(\bfu_\iota + \bfx) = \bfu_\iota + 2\bfx$ for all $\bfx \in \mathbb{R}^3$. Let also $P_\iota = \tau_\iota^{-1}(P)$. So each $P_\iota$ is a homothetic copy of $P$, scaled down by a factor of $1/2$, with $\bfu_\iota$ being the centre of the homothety. Over all vertices $\bfu_\iota$ of $P$, we get eight such homothetic copies $P_\iota$. Together, they form a partitioning of $P$, in the sense that they cover $P$ and they only meet at their boundaries.

\begin{lemma} \label{ping} Suppose that there is an integer point in the interior of $P$ outside of its diagonal $\bfo\bfu_{123}$, and also that there are no integer points on the boundary of $P$ other than its vertices. Then there is an integer point in the interior of one among $P_1$, $P_2$, $P_3$. \end{lemma}

\begin{myproof} Let $\bfw$ be the integer point in the interior of $P$ furthest away from the line through $\bfo$ and $\bfu_{123}$. Say it lies in $P_\iota$. Observe that $\bfw$ must lie in the interior of $P_\iota$, as otherwise $\tau_\iota(\bfw)$ becomes an integer point on the boundary of $P$ distinct from its vertices.

Suppose, for the sake of contradiction, that $\iota \in \{\varnothing, 123\}$. Then $\tau_\iota(\bfw)$ becomes an integer point in the interior of $P$ which lies twice as far away from the line through $\bfo$ and $\bfu_{123}$. This contradicts our choice of $\bfw$. So it must hold that $\iota \in \{1, 2, 3, 12, 13, 23\}$. We are done when $\iota \in \{1, 2, 3\}$. Otherwise, we reflect $\bfw$ across the centre of $P$. \end{myproof}

\begin{lemma} \label{bagel} Suppose that there is an integer point in the interior of one among $P_1$, $P_2$, $P_3$. Then there is also an integer point in $H$ distinct from $\bfo$, $\pm \bfu_1$, $\pm \bfu_2$, $\pm \bfu_3$. \end{lemma}

\begin{myproof} Suppose, for concreteness, that $\bfw = \langle 1 - w_1, w_2, w_3 \rangle$ is an integer point in the interior of $P_1$. Then $0 < w_i < 1/2$ for all $i$.

We proceed now to construct an integer point $\bfw' = \langle 1 - w'_1, w'_2, w'_3 \rangle$ in the tetrahedron $\bfo\bfu_1\bfu_{12}\bfu_{13}$, distinct from its vertices. This tetrahedron is the convex hull of $\bfu_1$ and its three neighbouring vertices of $P$. Or, from the perspective of $\Lambda$, it is the set of all points $\langle 1 - \lambda_1, \lambda_2, \lambda_3 \rangle$ with the $\lambda$'s nonnegative and $\lambda_1 + \lambda_2 + \lambda_3 \le 1$.

If $w_1 + w_2 + w_3 \le 1$, we set $\bfw' = \bfw$. Otherwise, we apply $\tau_1$ to $\bfw$ and then reflect the resulting image across the centre of $P$. This yields an integer point $\bfw'$ with $w'_i = 1 - 2w_i > 0$ for all $i$, as well as with $w'_1 + w'_2 + w'_3 = 3 - 2(w_1 + w_2 + w_3) < 1$.

Finally, let $\bfw'' = \bfw' - \bfu_1$. Then $\bfw''$ is an integer point in the tetrahedron $\bfo\bfu'_1\bfu_2\bfu_3$, distinct from its vertices. (Recall that $\bfu'_1$ is our shorthand notation for $-\bfu_1$.) However, this is one of the six constituent tetrahedra of $H$. \end{myproof}

Together, Lemmas \ref{ping} and \ref{bagel} complete the resolution of the non-diagonal case.

\section{Three-Dimensional Grids III} \label{3diii}

Suppose, throughout this section, that all integer points in the interior of $P$ lie on its diagonal $\bfo\bfu_{123}$, and also that there are no integer points on the boundary of $P$ other than its vertices.

We begin with some general remarks. Consider any parallelepiped $\Pi(\bfv_1, \bfv_2, \bfv_3)$ with integer and linearly independent $\bfv_1$, $\bfv_2$, $\bfv_3$. Let $I_0$ be the set of all integer points in its interior, $I_1$ the set of all integer points in the interiors of its faces, $I_2$ the set of all integer points in the interiors of its edges, and $I_3$ the set of its vertices. Let also $\Gamma$ be the lattice generated by $\bfv_1$, $\bfv_2$, $\bfv_3$.

The translation copies of $\Pi(\bfv_1, \bfv_2, \bfv_3)$ by the elements of $\Gamma$ form a tiling of $\mathbb{R}^3$. The number of tiles that each integer point belongs to is determined by its congruence class modulo $\Gamma$. This congruence class meets exactly one among $I_0$, $I_1$, $I_2$, $I_3$, and when it meets $I_i$ the number of tiles containing that integer point will be $2^i$. By estimating the number of integer points in a large ball, and letting the radius of that ball tend to infinity, we arrive at $\vol \Pi(\bfv_1, \bfv_2, \bfv_3) = \sum 1/2^i \cdot |I_i|$.

We now specialise this analysis to $P = \Pi(\bfu_1, \bfu_2, \bfu_3)$. Then $I_1 = I_2 = \varnothing$ and, as with every parallelepiped, $|I_3| = 8$. So there are exactly $p - 1$ integer points in the interior of $P$.

Since all of these points lie on the diagonal $\bfo\bfu_{123}$, we can write them down as $\mu_1\bfu_{123}$, $\mu_2\bfu_{123}$, $\ldots$, $\mu_{p - 1}\bfu_{123}$ with $0 < \mu_1 < \mu_2 < \cdots < \mu_{p - 1} < 1$. Observe, however, that the sum of any two integer points in $P$, taken modulo $\Gamma$, is an integer point in $P$ once again. So the set $\{0, \mu_1, \mu_2, \ldots, \mu_{p - 1}\}$ is closed with respect to addition modulo unity. Or, equivalently, the integer translates of $0$, $\mu_1$, $\mu_2$, $\ldots$, $\mu_{p - 1}$ form a lattice in $\mathbb{R}$. We conclude that $\mu_i = i/p$ for all $i$.

Since $-\bfu_4 = \langle s_1/s_4, s_2/s_4, s_3/s_4 \rangle$ belongs to $P$, we get that $\bfu_4 = -i/p \cdot \bfu_{123}$ for some $i$ with $1 \le i \le p$. So $-1/p \cdot \bfu_{123}$ lies on the segment $\bfo\bfu_4$, and hence it belongs to $B$. We can assume without loss of generality that in fact $\bfu_4 = -1/p \cdot \bfu_{123}$. (Otherwise, we can shrink our enveloping quadruple by replacing $\bfu_4$ with $-1/p \cdot \bfu_{123}$.) Thus $\bfu_1 + \bfu_2 + \bfu_3 + p\bfu_4 = \bfo$.

\begin{lemma} \label{1234} Each one of $\bfu_1 + \bfu_4$, $\bfu_2 + \bfu_4$, $\bfu_3 + \bfu_4$ is excellent. \end{lemma}

\begin{myproof} Consider $\bfv = \bfu_1 + \bfu_4$, for concreteness. Then $-\bfv$ is good because its replacement set is the triple $\{\bfu_1, \bfu_4, -\bfv\}$. For $\bfv$, the interesting case is when its replacement set is a quadruple. Then its replacement figure must be one of the two tetrahedra $\bfu_1\bfu_2\bfu_3\bfv$ and $\bfu_2\bfu_3\bfu_4\bfv$. By the linearity of the signed volume, $\vol \bfu_1\bfu_2\bfu_3\bfo + \vol \bfu_1\bfu_2\bfu_3\bfv = \vol \bfu_1\bfu_2\bfu_3\bfu_1 + \vol \bfu_1\bfu_2\bfu_3\bfu_4$, and so $\vol \bfu_1\bfu_2\bfu_3\bfv = t - s_4 < t$. Similarly, $\vol \bfu_2\bfu_3\bfu_4\bfv = t - s_1 < t$ as well. \end{myproof}

\begin{lemma} \label{sign} Suppose that $\bfu_1 + \bfu_2 + \bfu_3 + p\bfu_4 = \bfo$. Suppose also that neither one of $\bfu_1 + \bfu_4$, $\bfu_2 + \bfu_4$, $\bfu_3 + \bfu_4$ belongs to $B$. Then, without loss of generality, the signatures of $\bfu_1$, $\bfu_2$, $\bfu_3$, $\bfu_4$ are ${-}{+}{+}$, ${+}{-}{+}$, ${+}{+}{-}$, ${-}{-}{-}$, respectively. \end{lemma}

\begin{myproof} We say that coordinate position $i$ \emph{breaks} the sum $\bfu_j + \bfu_4$ when the absolute value of the $i$-th component of that sum is too large for it to belong to $B$.

We claim that each coordinate position breaks at most one sum, and furthermore that breakage can occur only when the $i$-th components of $\bfu_j$ and $\bfu_4$ are of the same nonzero sign while the $i$-th components of the other two $\bfu$'s are of the opposite nonzero sign. Clearly, this suffices: We get three coordinate positions to break three sums, and for each one of these three coordinate positions we can flip all signs at it without loss of generality.

Fix the coordinate position $i$ under consideration, and let $u_1$, $u_2$, $u_3$, $u_4$ be the $i$-th components of $\bfu_1$, $\bfu_2$, $\bfu_3$, $\bfu_4$, respectively. Suppose, for concreteness, that it is the sum $\bfu_1 + \bfu_4$ which is being broken. Then $u_1$ and $u_4$ must be of the same nonzero sign, as otherwise $|u_1 + u_4| \le \max\{|u_1|, |u_4|\}$. Suppose, for concreteness, that both of them are positive.

Consider $u_2$ and $u_3$ next. Suppose, for the sake of contradiction, that $u_2$ is nonnegative. Since now all three of $u_1$, $u_2$, $u_4$ are nonnegative and $u_1 + u_2 + u_3 + pu_4 = 0$, we get that $u_1 + u_4 \le u_1 + u_2 + pu_4 = -u_3$ and $|u_1 + u_4| \le |u_3|$. But this is a contradiction with the breakage of $\bfu_1 + \bfu_4$ at the coordinate position $i$ under consideration. So $u_2$ must be negative. Similarly, $u_3$ must be negative as well. \end{myproof}

\begin{lemma} \label{vol} Let $\bfv_1$, $\bfv_2$, $\bfv_3$ be three integer points with signatures ${-}{+}{+}$, ${+}{-}{+}$, ${+}{+}{-}$ whose sum $(a, b, c) = \bfv_1 + \bfv_2 + \bfv_3$ is of signature ${+}{+}{+}$. Then \[\vol \Pi(\bfv_1, \bfv_2, \bfv_3) \ge \min\{(a + 1)(b + c), (b + 1)(c + a), (c + 1)(a + b)\}.\] \end{lemma}

\begin{myproof} Let $\bfv_1 = (-a_1, b_1, c_1)$, $\bfv_2 = (a_2, -b_2, c_2)$, $\bfv_3 = (a_3, b_3, -c_3)$. Consider the multilinear optimisation problem with $9$ variables $a_1$, $b_1$, $\ldots$, $c_3$ and $12$ constraints $a_1 \ge 1$, $b_1 \ge 1$, $\ldots$, $c_3 \ge 1$ as well as $-a_1 + a_2 + a_3 = a$, $b_1 - b_2 + b_3 = b$, $c_1 + c_2 - c_3 = c$ where we wish to minimise the quantity \[\vs \Pi(\bfv_1, \bfv_2, \bfv_3) = \det \left(\begin{smallmatrix} -a_1 & b_1 & c_1\\ a_2 & -b_2 & c_2\\ a_3 & b_3 & -c_3 \end{smallmatrix}\right).\]

Suppose first that we replace $a_2$ and $a_3$ with $a_2 + r$ and $a_3 - r$, where $r \in [-(a_2 - 1); a_3 - 1]$. Since the altered $\vs \Pi(\bfv_1, \bfv_2, \bfv_3)$ is a linear function of $r$, it attains its minimum when $r$ is at one of its two extreme values. So, without loss of generality, either $a_2 = 1$ or $a_3 = 1$.

The two cases behave analogously, and so for the next step of the argument we focus on the scenario where $a_3 = 1$. Suppose now that we replace $a_1$ and $a_2$ with $a_1 + r$ and $a_2 + r$, where $r \ge -(a_1 - 1)$. The altered $\vs \Pi(\bfv_1, \bfv_2, \bfv_3)$ becomes a linear function of $r$ where the coefficient of $r$ equals \begin{gather*} -\det \left(\begin{smallmatrix} -b_2 & c_2\\ b_3 & -c_3 \end{smallmatrix}\right) + \det \left(\begin{smallmatrix} b_3 & -c_3\\ b_1 & c_1 \end{smallmatrix}\right) = \det \left(\begin{smallmatrix} b_3 & -c_3\\ -b_2 & c_2 \end{smallmatrix}\right) + \det \left(\begin{smallmatrix} b_3 & -c_3\\ b_1 & c_1 \end{smallmatrix}\right) = {}\\ {} = \det \left(\begin{smallmatrix} b_3 & -c_3\\ b_1 - b_2 & c_1 + c_2 \end{smallmatrix}\right) = \det \left(\begin{smallmatrix} b_3 & -c_3\\ b_1 - b_2 + b_3 & c_1 + c_2 - c_3 \end{smallmatrix}\right) = \det \left(\begin{smallmatrix} b_3 & -c_3\\ b & c \end{smallmatrix}\right). \end{gather*}

The last expression is clearly positive. So the altered $\vs \Pi(\bfv_1, \bfv_2, \bfv_3)$ is minimised when $r$ attains its smallest value; i.e., when $a_1 = 1$ and $a_2 = a$. Similarly, in the scenario where $a_2 = 1$, we can assume without loss of generality that $a_1 = 1$ and $a_3 = a$.

The same reasoning applies to $b_1$, $b_2$, $b_3$ and $c_1$, $c_2$, $c_3$. Hence, without loss of generality, $a_1 = b_2 = c_3 = 1$ as well as $\{a_2, a_3\} = \{1, a\}$, $\{b_1, b_3\} = \{1, b\}$, $\{c_1, c_2\} = \{1, c\}$. There are now just eight different options to sort through. The corresponding determinants work out to a total of four distinct values: $(a + 1)(b + c)$, $(b + 1)(c + a)$, $(c + 1)(a + b)$, $abc + a + b + c$. However, $(a + 1)(b + c) \le abc + a + b + c$, as the inequality simplifies to $0 \le a(b - 1)(c - 1)$. Thus the lattermost value can be safely discarded. \end{myproof}

We are only left to put the pieces together. By Lemma \ref{1234}, we are done when one of $\bfu_1 + \bfu_4$, $\bfu_2 + \bfu_4$, $\bfu_3 + \bfu_4$ belongs to $B$. Suppose, for the sake of contradiction, that neither one of them does. Then, without loss of generality, the signatures of $\bfu_1$, $\bfu_2$, $\bfu_3$, $\bfu_4$ must be as in Lemma \ref{sign}.

Let $\bfu_{123} = (a, b, c)$. Since $\bfu_4 = -1/p \cdot \bfu_{123}$, we get that $p$ divides all three of $a$, $b$, $c$. By Lemma \ref{vol}, it follows that $\vol P \ge 2p(p + 1)$. However, $p = \vol P$ by definition. This is our desired contradiction, as clearly $p < 2p(p + 1)$. The diagonal case has been resolved, too, and the proof of Theorem \ref{3dt} is complete.

\section{Higher Dimensions and Further Work} \label{further}

We conclude with a look into higher-dimensional grids.

First we describe one unweighted translation order on the $2 \times 2 \times 2 \times 3$ grid, in terms of the black-and-white colouring of $D = [-1; 1]^3 \times [-2; 2]$. Let $f(\bfx) = x_1 + 3x_2 + 6x_3 + 2x_4$. We colour all points $\bfx$ of $D$ with $f(\bfx) < 0$ in black, and all with $f(\bfx) > 0$ in white. This leaves eight nonzero points of $D$ uncoloured. These points split uniquely into two quadruples with sum $\bfo$ symmetric across the origin. We make one quadruple white, and the other one black.

This is enough to show that Theorem \ref{3dt} does not generalise to four or more dimensions. Notice, though, that the black and white points of $D$ in the counterexample can still be separated nonstrictly. This is somewhat unsatisfactory, as one might argue that the induced order is ``almost'' weighted. We call translation orders where both strict and nonstrict separation are impossible \emph{strongly unweighted}.

We proceed to exhibit one such translation order on the grid of size $3 \times 3 \times 3 \times 3$, with $D = [-2; 2]^4$. Consider $f(\bfx) = 4x_1 + 14x_2 + 21x_3 + 61x_4$. Direct computation shows that, over $D$, the only point where $f$ vanishes is the origin, while the values $\pm 2$ are not attained at all. We colour all points $\bfx$ of $D$ with $f(\bfx) \le -3$ or $f(\bfx) = 1$ in black, and all with $f(\bfx) = -1$ or $f(\bfx) \ge 3$ in white.

Is it true that, for all dimensions $d \ge 4$ and all positive integers $n$, there exists a strongly unweighted translation order on some $d$-dimensional grid all of whose dimensions exceed $n$? This does seem plausible. However, explicit constructions which would demonstrate it in full generality are not immediately obvious.

What can be said about the structure of the unweighted, or strongly unweighted, translation orders on higher-dimensional grids? For example, do they always admit simple descriptions similar to the two we just presented? Of course, a lot is going to depend on how exactly the notion of simplicity has been formalised.

Finally, there is also the even more open-ended question of finding other natural and interesting families of integer point sets $S$ for which all translation orders are weighted.

\section*{Acknowledgements}

The author is thankful to Ankan Bhattacharya for telling him about the original conjecture, as well as to Kevin Ren and Brice Huang for kindly agreeing to share their proof.

The present paper was written in the course of the author's PhD studies under the supervision of Professor Imre Leader. The author is also thankful to Prof.\ Leader for his unwavering support.

\end{document}